# Identification of Structural Model for Chaotic Systems

## Evgeny V. Nikulchev, Oleg V. Kozlov

Moscow Technological Institute, Moscow, Russia
Email: nikulchev@mail.ru

## ABSTRACT

This article is talking about the study constructive method of structural identification systems with chaotic dynamics. It is shown that the reconstructed attractors are a source of information not only about the dynamics but also on the basis of the attractors which can be identified and the mere sight of models. It is known that the knowledge of the symmetry group allows you to specify the form of a minimal system. Forming a group transformation can be found in the reconstructed attractor. The affine system as the basic model is selected. Type of a nonlinear system is the subject of calculations. A theoretical analysis is performed and proof of the possibility of constructing models in the central invariant manifold reduced. This developed algorithm for determining the observed symmetry in the attractor. The results of identification used in real systems are an application.

**Keywords:** Chaotic Dynamics; Attractors; Symmetry; Identification

## 1. Introduction

The first work on the reconstruction of the strange attractor from the time series has been publishing the results on hydrodynamics [1]. The article shows that you can get a satisfactory picture of the strange attractor of the geometric dimensions of a small, if the variables $x$, appearing in the equations of the dynamical system $dx/dt = F(x)$, use the m-dimensional vectors, derived from the elements of time series of the same principle, which in the problems of autoregression. That same year, F. Takens reported on his theorem, which was published a year later [2]. That it is the basis of all algorithms for time series analysis methods of nonlinear dynamics. The problem of determining the form of a dynamical system from its one-dimensional realization belongs to a class of incorrect problems. Unlike the problem of analyzing this issue is ambiguous, since there are infinitely many dynamical systems of various kinds which can play the existing signal with a given degree of accuracy.

The method of the global reconstruction of a dynamic system of equations for its one-dimensional realization was proposed in [3,4]. The algorithm is as follows. One-dimensional realization of the process is in a system, which is considered a "black box" recovered phase portrait on the Takens theorem, topologically equivalent to the attractor of the original system. A priori given equation is the method of least squares which is a set of un-known coefficients.

Now there are considerable publications, developing and constantly improving the proposed the method [5-11]. For example, in R. Brown and others [12] to reconstruction dynamic equations on the experimental time series with a broadband continuous spectrum use additional information about the dynamic and statistical properties of the original system contained in the implementation. In obtaining equation takes into account the values of Lyapunov exponents and the probability density, calculated from the original time series. However, the resulting evolution equations have a very cumbersome, inconvenient to use. In [13] used the hidden variables to write a model equation. In [12] describes a method for synchronizing the model with the original data. In a number of O.L. Anosov proposed reconstruction algorithm scalar differential equations for systems with delay.

However, the feature of many studies is that the proposed methods are illustrated with examples of simple low-dimensional model systems when we know in advance what should be the result of global reconstruction. It does not show substantial benefits given by the modifications of the basic method [3]. Described in the publications of the algorithms are tested on a number of known model systems with small dimension and a simple form right sides. The efficiency of the method is demonstrated by the example of time series generated by the real "black boxes".

The challenge is the need to work with noisy data in the processing of experimental time series. On the one hand, more desirable is the use of sequential differentiation to restore the phase trajectory, because it can get a model that contains, in general, approximately $n$ times smaller than the coefficients of the various non-linear than when using the method of delays. But differentiation will inevitably lead to increased noise components of high order. Without pre-filtering the time dependence of the second derivative can be noise-like process. In addition, methods of attachment are obvious flaws in the analysis significantly heterogeneous implementations, *i.e.*, signals in which areas with fast motion alternating with areas of slow motions.

Arbitrary choice of nonlinearities, as a rule, does not allow for a successful reconstruction of the dynamic equations for real systems. In particular, in [14] indicated the presence of three typical cases:

1) Reconstruction of locally describe the phase trajectory of the original system. In this case, the reconstructed model is unstable in the sense that the solution of these equations reproduces the signal under investigation only in a short period of time.

2) There is poor local predictability of the phase trajectory, but there is visual similarity of the phase portraits. Reconstruction of solution is stable in the sense of Poisson. In this case, the attractor of the reconstructed model has metric characteristics similar to those of the original attractor.

3) There is a good local predictability of the phase trajectory of any point in time values exceeding the characteristic correlation time. The phase portrait reconstructed model is identical to the original, and the system is Poisson stable.

## 2. Group Properties for the Construction of Equivalent Mapping

Consider the research questions of the global structure of orbits of dynamical systems, which do not depend on the choice of coordinate system, *i.e.*, admit a symmetry transformation. From a global point of view, the change of coordinates is a diffeomorphism (in the case of a smooth structure) or a homeomorphism (in the topological situation) phase spaces. Thus, we can introduce a natural equivalence relation between the dynamical systems associated with different classes of coordinate changes, and interpret the problem of describing the structure of orbits as a classification problem of dynamical systems up to these equivalence relations.

Given the equivalence of the trajectories on the torus [15], we can show that every periodic point of the map, whose spectrum does not contain a unit, defines certain $\mathbb{C}^1$-modules. Since these periodic orbits are separated from each other, their spectra can be outraged, regardless, at least, for any finite set of points. The modules from a variety of periodic orbits a independent.

On the other hand, at least in some cases, the local spectrum is a complete invariant with respect to a smooth conjugacy. Different approaches to the problem of local smooth conjugacy shown in [16].

For the global structure of orbits there are ways to build structures of independent modules associated with periodic orbits. In the case of infinitely many periodic orbits, like expanding mapping and hyperbolic automorphism of the torus, there are infinitely many invariants of the local $\mathbb{C}^1$-equivalence. For these two cases, which represent the simplest examples of hyperbolic systems, the spectra of periodic orbits form a complete system of invariants for $\mathbb{C}^1$- and $\mathbb{C}^\infty$-equivalences in the vicinity of these maps, respectively. Satisfactory description of the set of possible values for the eigen values of periodic points remains an open problem.

A variety of modules offer a significant, although not complete information about the smooth equivalence in the vicinity of the rotation. Every real analytic compression stores a unique way an affine structure. For the discrete map $\psi$, two structures identified near the ends of the segment are found in the middle. The transition functions between the two structures in any of the fundamental domain $\left[a, \psi(a)\right]$ generate infinite moduli space $\psi$ [16]. In practice, this can be interpreted to mean that there are changes of coordinates that will $\psi$ affine mapping from $\left[0, \psi^{-1}(0)\right]$ in the $[0, a]$ and $\left[\psi(a),1\right]$ in the $\left[\psi^2(a),1\right]$. Coordinates are uniquely determined up to two factors, one at each end. Then the map $\left[\psi(a),1\right] \rightarrow \left[\psi^2(a),1\right]$ can be normalized so that $\psi_a : [0,1] \rightarrow [0,1]$:

$$\psi_a = \frac{\psi\left(a + t\left(\psi(a) - a\right)\right) - \psi(a)}{\psi^2(a) - \psi(a)},$$

which can be extended to the whole real line by the formula

$$\psi_a\left(T + k\right) = \psi_a\left(T\right) + k, \ \ k \in \mathbb{Z}.$$

Thus, the two diffeomorphism $\psi_1$ and $\psi_2$, for which

$$\psi_1(0) = \psi_2(0) = 0$$

are equivalent if

$$\psi_2(T) = \psi_1(T + s) - \psi_1(s)$$

for some s in [0,1]. Constructed above the map transition depends on the factors that determine the linearization $\psi$ in the vicinity of the endpoints, as well as the choice of the reference point $a$ as well. It is obvious that the change does not alter the linearization map transition if $a$ amended accordingly. Changing the reference point leads to the substitution mapping the transition to the equiva-

　　　　　　　　　　　　　　　　　　　　　　　　　　　　　　　　　　　　　　*JMP*

lent.

The notion of smooth equivalence in accordance with [16] can be easily moved in case of continuous time. Flow equivalence is a contingency flows as a differentiable actions of a group $\mathbb{R}$ of real numbers. The structure of the orbits of the flow, unlike the case of discrete-time systems, has two distinct aspects: 1) The relative behavior of points on different orbits and 2) evolution of the initial condition along the orbit over time. There is a natural way to change the thread, keeping the first aspect of the structure of its orbits, and it does not change its orbit.

For the two streams $\varphi^t$ and $\varepsilon^t$ will define the time change by using the vector fields of infinitesimal generator of the symmetry group:

$$\xi = \frac{d\varphi^t}{dt}\bigg|_{t=0} \quad \text{and} \quad \eta = \frac{d\varepsilon^t}{dt}\bigg|_{t=0} \,.$$

From the uniqueness of solutions of differential equations, it follows that the zeros of the vector fields are fixed points of the corresponding flow. Thus, we conclude that $\xi(x) = 0$ if and only if $\eta(x) = 0$. In addition, if $x$ is not a fixed point, then the tangent vector to the curve $\{\varphi^t(x)\}$ and $\{\varepsilon^t(x)\}$, does not vanish and have the same direction.

It is natural to try to describe all the changes of time given flow modulo trivial substitutions. This problem is essentially equivalent to the problem of describing the space of all sufficiently smooth positive functions up to adding features that are derivatives of other smooth functions in the direction of flow.

It makes sense when determining equivalence to consider the conservation of structural stability for flows. One way of determining is based on the equivalence of all the perturbations [16]. Not being fully degenerate, this requirement is rarely performed, for example, in the presence of periodic orbits of periods are modules in this sense. We assume that the local topological equivalence should preserve the structural stability, thus transforming the homeomorphism can be quite close to the identity for small perturbations.

For all the above concepts in this section of the compact phase spaces is immaterial. In addition, a natural way to modify these definitions for the cases where for some points of the dynamical system is defined only for a finite time interval. Such a generalization leads to the concepts of local and semilocal structural stability.

For any contraction phase space can not have a smooth structure, so that these concepts are not directly applicable. However, the contraction mapping of a small disc in the Euclidean space, structural stability, as well as a hyperbolic linear map in the vicinity of the fixed point.

Symmetry rotation is not structurally stable. Since the topological conjugacy preserves the periodic orbits, and the transformation of rotation by an irrational angle can not be associated with the transformation of rotation, for which the corresponding number is rational. But since the rational numbers, and the set of irrational numbers are dense, then some arbitrarily small perturbations of the turn to rational angle you can find the rotation by an irrational angle, and vice versa.

## 3. Geometric Method for the Reduction of Systems

According to [15] on centrally-stable manifold, the system around an equilibrium point $O$ takes the form

$$\begin{aligned}
\dot{y} &= By + f_1(x, y, z), \\
\dot{z} &= Cz + f_2(x, y, z), \\
\dot{x} &= Ax + \psi_0(x, y, z),
\end{aligned} \tag{1}$$

where

$$x \in \mathbb{R}^m, \ y \in \mathbb{R}^k, \ z \in \mathbb{R}^{n-m-k} \quad \left(j = \overline{1, m}\right);$$

$$\text{spectr } A = \{\lambda_{m+1}, \cdots, \lambda_k\}, \quad \operatorname{Re} \lambda_j < 0, \quad \left(j = \overline{m+1, k}\right);$$

$$\text{spectr } C = \{\lambda_{k+1}, \cdots, \lambda_n\}, \quad \operatorname{Re} \lambda_j > 0, \quad \left(j = \overline{k+1, n}\right);$$

$\mathbb{C}^r$-functions $f_1, f_2$ and $\psi_0$ together with its first derivatives vanish at the origin.

In this case, the right side of the system may depend on the control u either continuously (in this case, smooth manifolds, discussed below, depend continuously on $u$), or smooth. In the latter case, u is included in the number of "central" variables $x$ and, thus, further investigated the variety and the foliation will have the smoothness of $u$, equal to the smoothness of $x$.

**Theorem 1.** In a small neighborhood of the equilibrium state there exists a $(m + k)$-dimensional center-stable invariant manifold. Central-stable manifold is not uniquely defined, but for any two manifolds $W_1^{sC}$ and $W_2^{sC}$ functions $\phi_1^{sC}$ and $\phi_2^{sC}$ the class $\mathbb{C}^r$ that contains the point $O$ and concerns at this point subspace $\{z = 0\}$, define the same the same symmetries group of the point $O$.

*Proof.* Without loss of generality, we consider the second-order system. Let the system defined on a manifold gives rise to a one-parameter transformation group generated the operator

$$A = \xi_1(x)\frac{\partial}{\partial x_1} + \xi_2(x)\frac{\partial}{\partial x_2},$$

with infinitesimal operator

$$X = \eta_1(x)\frac{\partial}{\partial x_1} + \xi_2(x)\frac{\partial}{\partial x_2}\,.$$

We need to show that the operator

$$\tilde{A} = \tilde{\xi}_1\left(x', \tau\right)\frac{\partial}{\partial x_1} + \tilde{\xi}_2\left(x', \tau\right)\frac{\partial}{\partial x_2}.$$

system defined by $W_2^{sC}$ an equivalent ($\tau$ - group option).

We write the transformation of the group as a Li series (the operator exponential):

$$x_1' = e^{\tau X} x_1,$$
$$x_2' = e^{\tau X} x_2,$$
$$x_1 = e^{\tau X} x_1',$$
$$x_2 = e^{\tau X} x_2'.$$

We represent $\tilde{A}$ in the previous coordinates. For this we compute:

$$A = \left(\tilde{A}x_1\right)\Big|_{\bar{x}=x'(x)}\frac{\partial}{\partial x_1} + \left(\tilde{A}x_2\right)\Big|_{\bar{x}=x'(x)}\frac{\partial}{\partial x_2}$$
$$= \left(\tilde{A}e^{-\tau X}x_1'\right)\Big|_{\bar{x}=x'(x)}\frac{\partial}{\partial x_1} + \left(\tilde{A}e^{-\tau X}x_2'\right)\Big|_{\bar{x}=x'(x)}\frac{\partial}{\partial x_2},$$

it follows that

$$\left(\tilde{A}e^{-\tau X}x_1'\right)\Big|_{\bar{x}=x'(x)} = \xi_1\left(x\right),$$
$$\left(\tilde{A}e^{-\tau X}x_2'\right)\Big|_{\bar{x}=x'(x)} = \xi_2\left(x\right).$$

Since $\xi_{1,2}$ do not depend on $\tau$, as determined by the operator $X$, we have

$$\frac{d}{d\tau}\left(\tilde{A}e^{-\tau X}x_1'\right) = 0,$$

so

$$\frac{\partial\tilde{A}}{\partial\tau}e^{-\tau X}x_1' - \tilde{A}Xe^{-\tau X}x_1' + X\tilde{A}e^{-\tau X}x_1' = 0.$$

A similar formula holds for the second coordinate, which together form a differential equation:

$$\frac{\partial\tilde{A}}{\partial\tau} = \tilde{A}X - X\tilde{A} = \left[\tilde{A}, X\right]$$

with initial condition

$$\tilde{A}\left(x', \tau\right)\Big|_{\tau=0} = A\left(x'\right).$$

Solution of the resulting Cauchy problem can be obtained by expanding the operator $\tilde{A}\left(x', \tau\right)$ in a Taylor series in powers of $\tau$:

$$\tilde{A}\left(x', \tau\right) = \tilde{A}\left(x'\right) + \tau\frac{\partial\tilde{A}}{\partial\tau}\Big|_{\tau=0} + \frac{\tau^2}{2!}\frac{\partial^2\tilde{A}}{\partial\tau^2}\Big|_{\tau=0} + \cdots,$$

thus, we have:

$$\frac{\partial\tilde{A}}{\partial\tau}\Big|_{\tau=0} = \left[A, X\right].$$

Similarly,

$$\frac{\partial^2\tilde{A}}{\partial\tau^2}\Big|_{\tau=0} = \left[\left[A, X\right], X\right].$$

Finally, a number will look like in the Hausdorff form:

$$\tilde{A} = A + \tau\left[A, X\right] + \frac{\tau^2}{2!}\left[\left[A, X\right], X\right], \cdots.$$

By Theorem from [15] function $\phi_1^{sC}$ and $\phi_2^{sC}$ determine the same expansion in Taylor series, which suggests field commutations:

$$\left[A, X\right] = 0,$$

*i.e.*

$$\tilde{A} = A.$$

Thus, the theorem is proved.

At the reference time $t \to -t$ of $A$, $B$ and $C$ go into respectively $A$, $B$ and $C$. Thus, part of the spectrum of characteristic exponents corresponding to the variable z, is now left of the imaginary axis, and part of the spectrum, corresponding to the variable $y$—the right of it. To the system obtained from system (1) by time reversal, we can again apply the theorem about the Central stable manifold and get the following theorem on center-unstable manifold.

**Theorem 2.** In a small neighborhood of the equilibrium state $O$ there exists an $(n - k)$-dimensional invariant manifold is a $\mathbb{C}^r$-smooth. Center unstable manifold includes all trajectories that remain in a small neighborhood of $O$ for all negative values of time. For any two $W_1^{usC}$ and $W_2^{usC}$ functions $\phi_1^{usC}$ and $\phi_2^{usC}$, the class $\mathbb{C}^r$ containing the point $O$ and tangent at this point subspace $\{y = 0\}$ define the same symmetries group.

Intersection of center-stable and unstable manifolds of a central $\mathbb{C}^r$-smooth m-dimensional invariant center manifold, defined by the equation of the form $(y, z) = \phi^C(x)$. Function $\phi^C$ with all derivatives, in particular, the Taylor expansion of functions $\phi^C$ at $O$ is uniquely determined by the system.

Straightening of the central stable and center-unstable manifolds, as well as rectification of the strong stable and strong unstable invariant foliations on these manifolds leads to the good fact.

**Theorem 3 [15].** With a $\mathbb{C}^{r-1}$-smooth transformation of the system (1) can locally be reduced to

$$\dot{y} = \left(A + F_1\left(x, y, z\right)\right)y,$$
$$\dot{z} = \left(C + F_2\left(x, y, z\right)\right)z,$$
$$\dot{x} = Bx + \Psi_0\left(x\right) + \Psi_1\left(x, y, z\right)y + \Psi_2\left(x, y, z\right)z,$$

where $\Psi_0$—is a $\mathbb{C}^r$-smooth function, which together with its first derivative vanishes at $x = 0$; $F_{1,2}$ are functions that vanish at the origin:

 

$$\Psi_{1,2} \in \mathbb{C}^{r-1}; \quad \Psi_1\left(x, y, 0\right) = \Psi_2\left(x, 0, z\right) = 0.$$

Here, the local center-unstable manifold given by the equation $\{y = 0\}$, the local center-stable manifold—$\{z = 0\}$, and a local center manifold—$\{y = 0, z = 0\}$. Strong stable foliation consists of surfaces $\{x = \text{const}, z = 0\}$, and the layers are strongly unstable foliation of the form $\{x = \text{const}, y = 0\}$.

A similar theory is constructed for discrete systems.

For systems that admit the group of symmetries and construct models of systems, reduced to a central invariant manifold results can be summarized as follows [17].

In the local region of the qualitative dynamic behavior of systems typologically equivalent system, reduced to the central manifold

$$\dot{x} = Ax\left(t\right) + \Psi_0\left(x\right), \quad (2)$$

where $\Psi_0\left(x\right)$—is a $\mathbb{C}^r$-smooth function

$$\left(\Psi_0\left(x_0\right) = \Psi_0'\left(x_0\right) = 0, x_0 = 0\right),$$

whose structure is determined by the symmetry transformation.

For a discrete system, the local center manifold is determined by the system:

$$x\left(t+1\right) = Ax\left(t\right) + \vartheta_0\left(x\right), \quad (3)$$

where $\vartheta_0\left(x\right)$—is a $\mathbb{C}^r$-smooth function

$$\left(\vartheta_0\left(x_0\right) = \vartheta_0'\left(x_0\right) = 0, x_0 = 0\right)$$

defined on the basis of symmetry groups constructed on the reconstructed attractor.

For the reconstruction of a nonlinear system in the forms (2), (3) proposed the allocation of local regions of phase trajectories $\varphi_1\left(x\right)$ and $\varphi_2\left(x\right)$ that are close to periodic, and the construction of finite-transformations transform one area to another. That is, the construction of the symmetry group of phase trajectories, which is characterized by the transformation of graphs:

$$\text{graph}\left\{\varphi_1\left(x\right)\right\} \rightarrow \text{graph}\left\{\varphi_2\left(x\right)\right\}.$$

The resulting transformations determine the structure of the desired evolution equations.

## 4. Algorithms and Applications

The result of reconstruction algorithms [1] is a set of points belonging to the attractor of the system. The task is an automatic search for symmetries of the local sections of the phase trajectories. We seek the symmetry of translation, rotation, stretching and compression.

You must select the initial sequence collection sites so that when reducing them to a single scale, position and angle of rotation, they would be most similar to each other, as well numerical indicators of transformations taking one piece to another without breaking of symmetry, to give a numerical estimate degree of symmetry breaking. The **Figure 1** shows a schematic illustration of the problem.

The task of allocating a set of similar fragments in the circuit is defined as follows:

*Input*: sequence of $n$ points $\left\langle v_1, v_2, \cdots, v_n \right\rangle$, the function $Dist\left(V_1, V_2\right)$, which gives an assessment of similarity between any two fragments of the original sequence

$$V_i = \left\langle v_j, v_{j+1}, \cdots, v_{j+p} \right\rangle.$$

*Output*: a set of disjoint fragments of the original sequence $S = \left(V_1, V_2, \cdots, V_k\right)$, such that for any other set of fragments $S' = \left(V_1, V_2, \cdots, V_h\right)$ of the condition:

$$\forall S' \neq S : \max_{V_i', V_j' \in S'} \left( \frac{1}{1 + Dist\left(V_i', V_j'\right)} \right) \cdot \sum_{i=1}^{k} |V_i'| \times \frac{1}{m}$$

$$\geq \max_{V_i, V_j \in S} \left( \frac{1}{1 + Dist\left(V_i, V_j\right)} \right) \cdot \sum_{i=1}^{h} |V_i| \cdot \frac{1}{m}$$

Number of different fragments, which can be identified in the original sequence:

$$m = \frac{1}{2} n\left(n-1\right)$$

Number of arbitrary sets of fragments, determined by the size of the set of all subsets of the set $\left\langle V_1, V_2, \cdots, V_m \right\rangle$:

$$u = P\left(m\right) = 2^m$$

Based on this challenge

$$T\left(n\right) = \Theta\left(2^{n^2}\right).$$

*Marking*. In the simplest case for the selection of fragments of trajectories satisfying the decision to consider the fragments located between all possible pairs of contour points. In most problems, this approach is not possible because of limited computer resources. So it makes sense to consider as a start / end fragment only some points of the contour, for example, extremes of functions of each of the coordinates and evenly distributed along the length of the markers on the slowly varying parts (**Figure 2**). This approach provides a sufficient rational partition of the circuit, although it suffers from some redundancy.

Partitioning algorithm can be anything important to exclude monotonous plots and identify the markers characteristic sections of the circuit. Extracting fragments. To isolate fragments of contour to study the simplest case, you can use brute force pairs of markers. In this case, the number of fragments needed to be addressed:

$$F_{\text{полн}} = \frac{n^2 - n}{2},$$

 

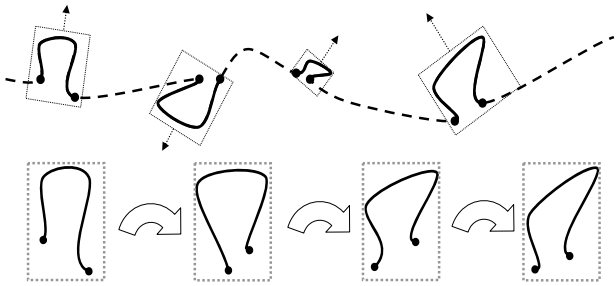

**Figure 1. Schematic illustration of the problem.**

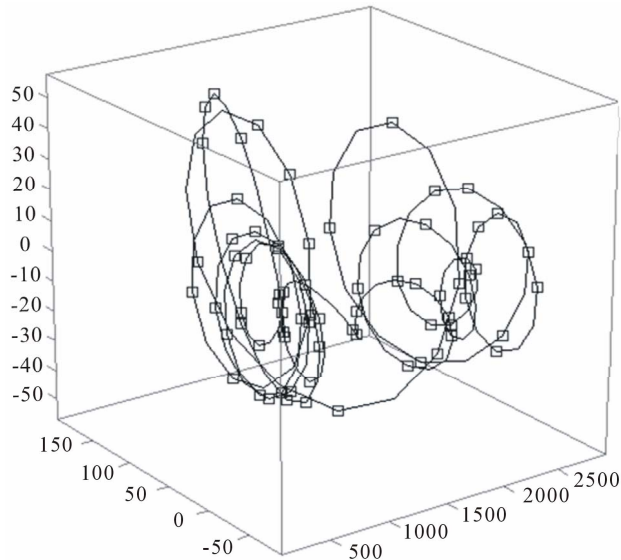

**Figure 2. Example of marking.**

where *n* — number of markers on the contour.

It makes sense to limit the length of a fragment in the original contour points as the top and bottom as well as fragments, the length of the approaching to the original contour, are not of interest to identify patterns of distribution symmetry, and too small fragments do not reflect the specific behavior of the system. In addition, the introduction of restrictions on the length of the above - increase the number of fragments for consideration is linear, as opposed to quadratic, with the full brute force.

Selected fragments for the possibility of a comparison between a need to interpolate in a way that they all described the same fixed number of points (**Figure 3**).

*Normalization*. After reducing to one the number of points each of the fragments should be subjected to the normalization procedure. The aim of the procedure is to convert the fragment into a descriptor—the image that is invariant with respect to the transport, rotation, and scaling the original fragment, and obtaining the numerical parameters of this transformation (**Figure 4**).

The proposed normalization procedure is a sequence of the following:

Isolation in a fragment of the so-called "axis"—pairs of points, the distance between the maximum in the current space projection.

Rotation of the fragment so that its axis were parallel to the coordinate axes.

The construction restricting fragment of the n-dimensional parallelepiped and the shift of the fragment so that the center of the parallelepiped coincides with the center coordinates.

Scaling the fragment so that the length of its main axis is equal to 1.

Fragment of an n-dimensional the contour *F* consisting of *m* real points, represented as a matrix $m \times n$:

$$F_b = \begin{bmatrix} X_1 & X_2 & \cdots & X_n \end{bmatrix},$$

where $X_i = \begin{bmatrix} x_{i,1}, x_{i,2}, \cdots, x_{i,m} \end{bmatrix}^T$.

For a compact description of the transformation fragment as a single matrix, we give it to homogeneous coordinates using matrix transformations, obtain:

$$F_c = \begin{bmatrix} X_1 & X_2 & \cdots & X_n & 1 \end{bmatrix}$$

In the normalization process uses three types of transformations: translation, scaling and rotation.

Translation transformation in a matrix form is:

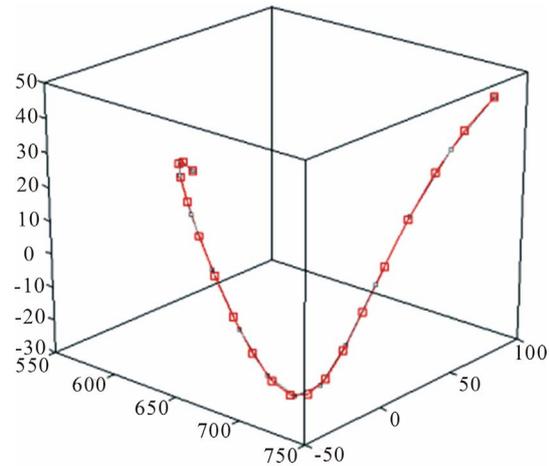

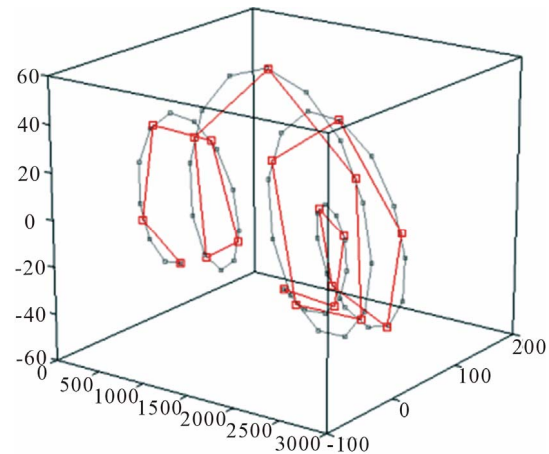

**Figure 3. Interpolation of selected fragments.**

                                                                          

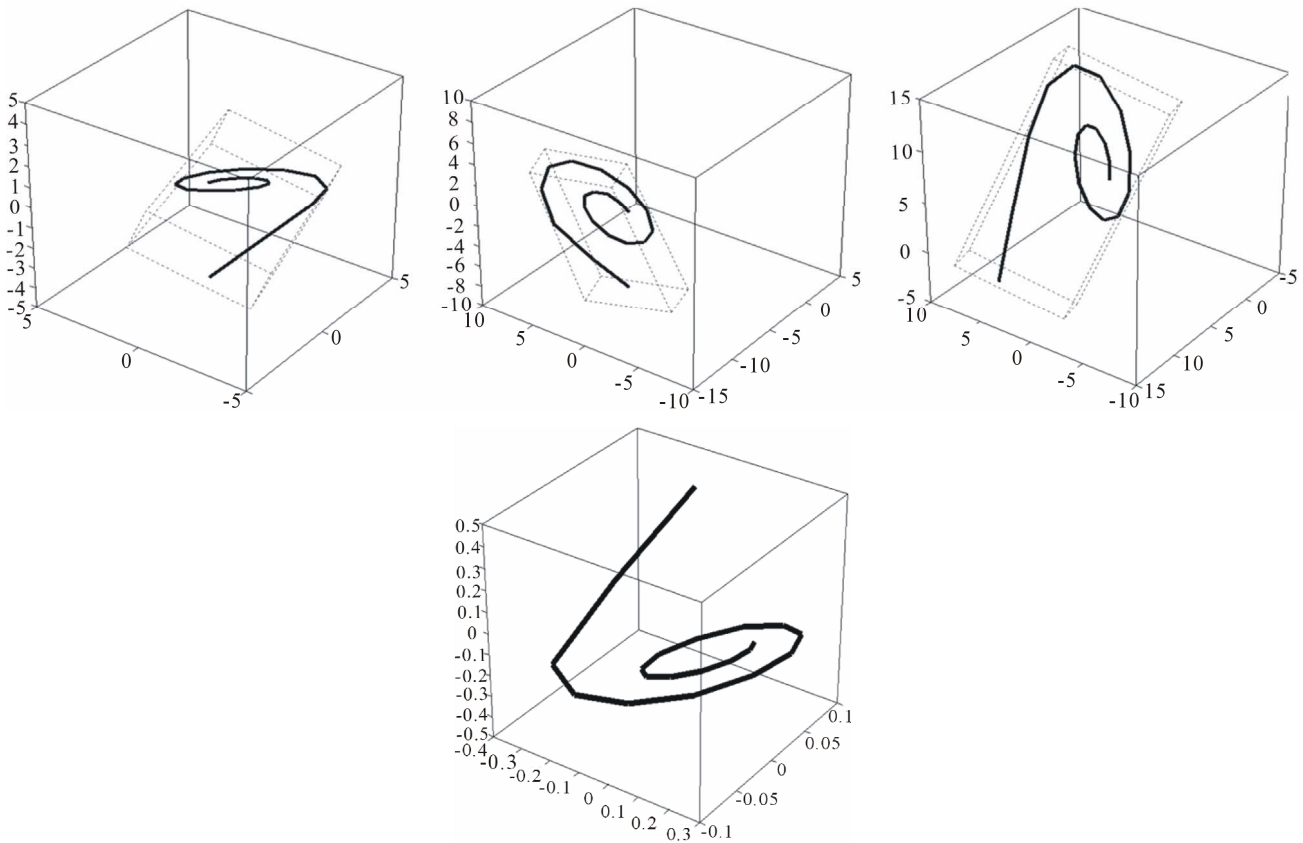

**Figure 4. Three initial symmetrical 3-dimensional fragment (top) and their total descriptor (bottom).**

$$M_{\text{sh}} = \begin{bmatrix} 1 & 0 & \cdots & 0 & 0 \\ 0 & 1 & \cdots & 0 & 0 \\ \vdots & \vdots & \ddots & \vdots & \vdots \\ 0 & 0 & \cdots & 1 & 0 \\ -g_{\text{sh},1} & -g_{\text{sh},2} & \cdots & -g_{\text{sh},n} & 1 \end{bmatrix}, \quad (4)$$

where $g_{\text{sh},i}$ $(i = 1, \cdots, n)$ —the shift along the axis $X_i$.

Transformation of the scaling down and scaling up

$$M_{\text{sc}} = \begin{bmatrix} g_{\text{sc}} & 0 & \cdots & 0 & 0 \\ 0 & g_{\text{sc}} & \cdots & 0 & 0 \\ \vdots & \vdots & \ddots & \vdots & \vdots \\ 0 & 0 & \cdots & g_{\text{sc}} & 0 \\ 0 & 0 & \cdots & 0 & 1 \end{bmatrix}, \quad (5)$$

where $g_{\text{sc}}$ —a scaling factor.

Rotation transformation is:

$$M_{\text{rt}} = \begin{bmatrix} 1 & 0 & \cdots & \cdots & 0 & 0 \\ 0 & \ddots & & & 0 & 0 \\ \vdots & & \cos(g_{\text{rt}}) & \sin(g_{\text{rt}}) & & \vdots \\ \vdots & & -\sin(g_{\text{rt}}) & \cos(g_{\text{rt}}) & & \vdots \\ 0 & 0 & & & \ddots & 0 \\ 0 & 0 & \cdots & \cdots & 0 & 1 \end{bmatrix}, \quad (6)$$

where $g_{\text{rt}}$ —angle, a significant elements of the matrix are located in cells $(i, i)$, $(i, j)$, $(j, i)$ и $(j, j)$.

## Normalization Algorithm Is as Follows

1) Initialize the total transformation matrix of normalization:

$$M_{\text{norm}} = \begin{bmatrix} 1 & 0 & \cdots & 0 \\ 0 & 1 & \cdots & 0 \\ \vdots & \vdots & \ddots & \vdots \\ 0 & 0 & \cdots & 1 \end{bmatrix}.$$

2) Setting the order of normalization $k = 1$.

3) $F^k = F_{\text{r}}$.

4) Finding the fragment $A^k$ —

$$A^k = \begin{bmatrix} x_{i,1} & x_{i,2} & \cdots & x_{i,n} & 1 \\ x_{i,1} & x_{j,2} & \cdots & x_{j,n} & 1 \end{bmatrix}$$

the largest segment connecting two points of the loop $F^k$ in space. With a simple enumeration pairs of points fragment of $F^k$ sought a pair with the greatest distance:

$$\max \sqrt{\sum_{d=k}^{n} (x_{d,j} - x_{d,i})^2}$$

5) Transfer $A^k$ so that it coincided with the first

                                                                                          

point of origin. Conversion is performed by multiplying with the matrix translation $M_{\text{sh},k}$ form (4) with coefficients

$$g_{\text{sh},i} = A^k\left(1,i\right):$$

$$A_{\text{sh}}^k = A^k \times M_{\text{sh},k}.$$

6) Complement the overall normalization of the current transformation matrix:

$$M_{\text{norm}} = M_{\text{norm}} \times M_{\text{sh},k}.$$

7) Rotation $A_{\text{sh}}^k$ so that it coincided with the direction of the axis $X_k$.

a) Initialize rotation matrix:

$$M_{\text{rt},k} = \begin{bmatrix} 1 & 0 & \cdots & 0 \\ 0 & 1 & \cdots & 0 \\ \vdots & \vdots & \ddots & \vdots \\ 0 & 0 & \cdots & 1 \end{bmatrix}$$

b) Setting the plane of rotation $p = n$.

$$A_{\text{rt},p}^k = A_{\text{sh}}^k$$

c) Calculating the angle between the projection $A_{\text{rt},p}^k$ on the plane $OX_1X_p$ and the axis $X_1$ of the formula:

$$\varphi_{k,p} = \begin{cases} \arctan\left(\dfrac{A^k(2,k)}{A^k(2,p)}\right), & A^k(2,k) > 0, A^k(2,p) > 0; \\[2mm] -\arctan\left(\dfrac{A^k(2,k)}{A^k(2,p)}\right) + \dfrac{\pi}{2}, & A^k(2,k) \le 0, A^k(2,p) > 0; \\[2mm] \arctan\left(\dfrac{A^k(2,k)}{A^k(2,p)}\right) + \pi, & A^k(2,k) \ge 0, A^k(2,p) \le 0; \\[2mm] -\arctan\left(\dfrac{A^k(2,k)}{A^k(2,p)}\right) + \dfrac{3\pi}{2}, & A^k(2,k) > 0, A^k(2,p) \le 0. \end{cases}$$

d) Rotation $A_{\text{rt},p}^k$ with a matrix form (6) by an angle

$$g_{\text{rt}} = -\varphi_{k,p}:$$

$$A_{\text{rt},p-1}^k = A_{\text{rt},p}^k \times M_{\text{rt},p}.$$

e) Complement the current rotation matrix transformation:

$$M_{\text{rt},k} = M_{\text{rt},k} \times M_{\text{rt},p}.$$

f) Setting the next plane of rotation $p = p - 1$.

g) If $p \ge 1$, then go to step {d}.

8) Complement the overall normalization of the current transformation matrix:

$$M_{\text{norm}} = M_{\text{norm}} \times M_{\text{rt},k}.$$

9) The transformation of the fragment by the current matrix of normalization:

$$F^{k+1} = F^k \times M_{\text{norm}}$$

10) Installation of the next order of normalization $k = k + 1$.

11) If $k < n$, then go to step {4}.

12) Formation of the bounding parallelepiped fragment $B$ is described by the matrix $2 \times n$:

$$B = \begin{bmatrix} b_{1,1} & b_{1,2} & \cdots & b_{1,n} \\ b_{2,1} & b_{2,2} & \cdots & b_{2,n} \end{bmatrix},$$

$$b_{1,i} = \min\left(X_i\right), b_{2,i} = \max\left(X_i\right).$$

13) Transfer fragment so that the center of the bounding parallelepiped it coincided with the origin. Conversion is performed by multiplying with the translation matrix $M_{\text{sh}}$ form (6) with coefficients

$$g_{\text{sh},i} = -\left(B(1,i) + B(2,i)\right)/2:$$

$$F_{\text{cn}} = F^k \times M_{\text{sc}}.$$

14) Scaling of the fragment in a way that restricts its parallelepiped fit into the unit cube. This conversion scaling matrix $M_{\text{sc}}$ form (5) with a coefficient

$$g_{\text{sc}} = B(2,i) - B(1,i):$$

$$F_{\text{norm}} = F_{\text{cn}} \times M_{\text{sc}}$$

15) Complement the overall normalization of the current transformation matrix:

$$M_{\text{norm}} = M_{\text{norm}} \times M_{\text{cn}} \times M_{\text{sc}}.$$

End.

The result of the normalization procedure is the descriptor the fragment $F_{\text{norm}}$, the matrix $M_{\text{norm}}$ is converted into a fragment of the descriptor and a set of indicators and all the intermediate transformation matrices: $n-1$ matrices

$$M_{\text{sh},k}, \ \left(n^2 - n\right)/2$$

matrices $M_{\text{rt},p}$ and one matrix $M_{\text{sc}}$.

In what is easily obtain the transformation matrix of one fragment to another:

$$M_{\text{trans}AB} = M_{\text{norm}A} \times M_{\text{norm}B}^{-1}.$$

*Evaluation of symmetry breaking.* After receiving the handles of the two supposedly symmetrical fragments $F_A$ and $F_B$ can give a numerical estimate of divergence between them $D_{AB}$, for it is proposed to use the properties of Fourier series. Also, this indicator can be used as an estimate of the degree of symmetry breakdown.

Successively translate the descriptor in both frequency representation, using the discrete Fourier transform for each point:

$$s_{d,k} = \sum_{p=1}^{m} x_{d,p} \mathrm{e}^{\frac{-2\pi i}{m}(k-1)(p-1)}, 1 \le d \le n.$$

For an $n$-dimensional the descriptor in this case applies $n$-dimensional Fourier transforms independently for each coordinate. Result of the transformation will range following structure:

$$S_d = \begin{bmatrix} s_{d,0} & s_{d,1} & \cdots & s_{d,m} \end{bmatrix}^T, 1 \le d \le n.$$

Couples $\left(S_1 \quad S_{n-1}\right), \left(S_2 \quad S_{n-2}\right), \cdots, \left(S_{(n-1)/2} \quad S_{n/2}\right)$ are complex conjugate numbers. The transition from the spectrum to the original contour by using the inverse discrete Fourier transform of the form:

$$x_{d,k} = \frac{1}{m} \sum_{p=1}^{m} s_{d,p} e^{\frac{2\pi i}{m}(k-1)(p-1)}, 1 \le d \le n.$$

For the purpose of comparing the contours of the Fourier transformation is convenient because it provides a selection of principal components. The closer a pair of conjugate elements in the middle of the spectrum, the smaller, high frequency, especially circuit it describes. Each of the pairs of conjugate elements of the spectrum determines when restoring an ellipse, and all recovered shape is a superposition of these ellipses. Each of the next pair will be added to restore the image all the more finer details of the original contour. The first, odd element of the spectrum indicates the position of the center of the figure for comparison; the contours are not being used.

To calculate the distance between the fragments is proposed to summarize the relevant elements of the divergence of the spectra of these contours, after Fourier transformation. Higher-frequency pair (close to the center of the spectrum), the elements of the spectrum have less impact on the indicator closeness, than the low frequency.

Corresponding criterion is:

$$D_{A,B} = \sum_{i=1}^{q} \left( \left( \sqrt{\sum_{j=1}^{n} I^2} + \sqrt{\sum_{j=1}^{n} R^2} \right) \cdot \beta_i \right)$$
$$I = \left( \operatorname{Im}\left(S_{B,i,j}\right) - \operatorname{Im}\left(S_{A,i,j}\right) \right),$$
$$R = \left( \operatorname{Re}\left(S_{B,i,j}\right) - \operatorname{Re}\left(S_{A,i,j}\right) \right)$$

where $q$—the number of conjugate pairs of elements of the spectrum, $n$—the dimension of phase space; $\beta_i$— discount factors $(i = 1, \cdots, q)$, which determine the degree of importance of components of the spectrum of frequencies—the relevant elements of the spectrum contours $A$ and $B$. Coefficients $\beta_i$ are chosen empirically, based on the conditions of the problem.

Once the fragments have been normalized circuit, must choose from a variety of sequences of one the most satisfying the requirements of the problem—containing the most balanced between fragments. Number sequences that can be composed of non-overlapping fragments, even for small problems with the number of markers about 100 very large, so it is proposed to solve by using a genetic algorithm.

## 5. Examples

### 5.1. Detection of Symmetries for a Model Example

Source data set consists of 300 points. These are the result of Ressler system simulation

$$\dot{x}_1 = -\left(x_2 + x_3\right),$$
$$\dot{x}_2 = x_1 + 0.2x_2,$$
$$\dot{x}_3 = 0.2 + x_3\left(x_1 - 2.6\right).$$

The observed parameter $x_3$ reconstructed attractor (**Figure 5**).

Attractor has been subjected to alignment markers. The total number of markers: 144. After labeling of the 10238 possible paths have been selected 1375, the others did not satisfy the restrictions in length: from 5 to 50 points of the original contour. Each segment was interpolated 60-th evenly spaced points along its length and subjected to the normalization procedure. After that was constructed adjacency matrix fragments (**Figure 6**). Been carried out selection of solutions by using genetic algorithm (population size 200, the probability of mutation rates $\alpha = 0.3$ and $\beta = 0.1$, limiting the number of steps lengths: 5, an elite selection of a new generation, the type of parental choice—outbreeding, unique individual in the population) and selected the winning solution.

The time of the preprocessing step: 3 min. 27 sec. Time step of genetic selection decisions: 4 min. 36 sec., the population converged 5 times for 533 iterations, the decision of the winner: the similarity of the fragments 0.3917, length 197.

The results are shown in **Figure 7**. Part of the trajectory for which were not symmetrical parts shown dotted.

### 5.2. Reconstruction of the Model to Experimental Data

Original data set consists of 200 points. The data represent the traffic network. Smoothing was not applied. The reconstructed attractor has been subjected to alignment markers. The total number of markers: 192. After labeling of the 18153 possible paths have been selected 2633, the others did not satisfy the restrictions in length: from 15 to 30 points of the original contour. Each segment was interpolated 60-th evenly spaced points along its length and subjected to the normalization procedure in **Figure 8**.

Was carried out the evolutionary selection of solutions (population size 200, the probability of mutation rates and limiting the number of steps lengths: 8, an elite se-lection of a new generation, the type of parental choice—outbreeding, unique individual in the population) and selected the winning solution. The time of the preprocessing stage: 12 min. 43 sec. Time step of genetic selection decisions: 6 min. 11 sec., The population converged 8 times for 652 iterations, the decision of the winner: the similarity of fragments of 0.576, length 125. The results are shown in **Figure 9**.

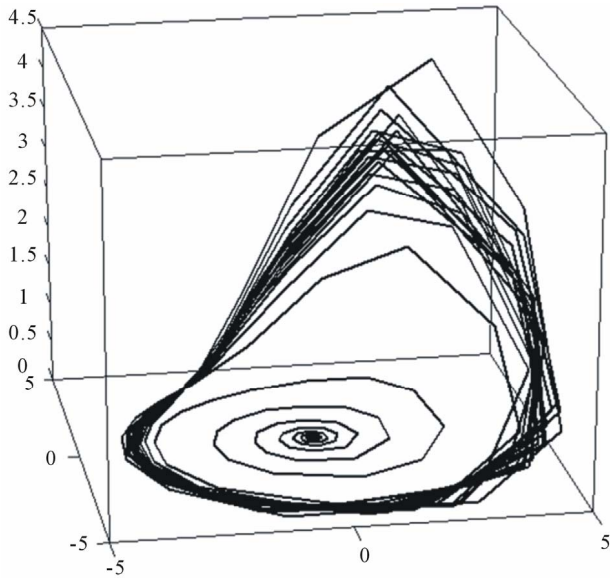

**Figure 5. Attractor.**

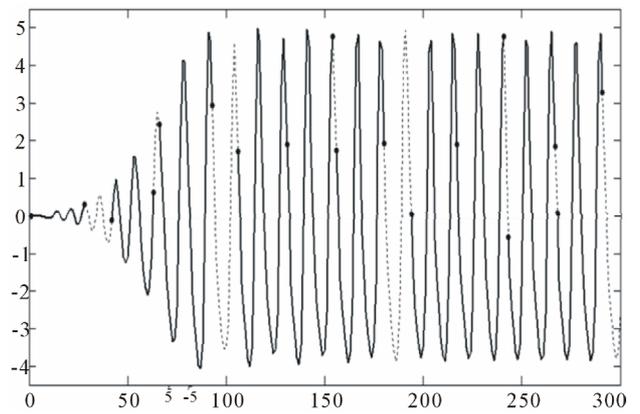

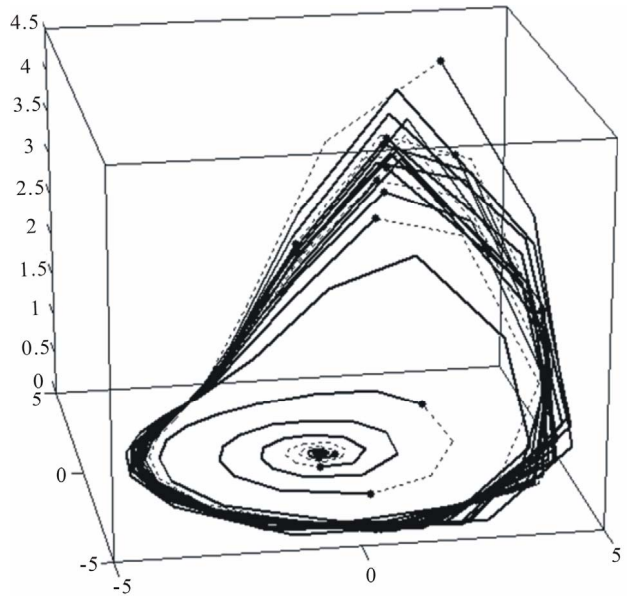

**Figure 7. Decision-winner in the time series (top) and in the attractor (botom).**

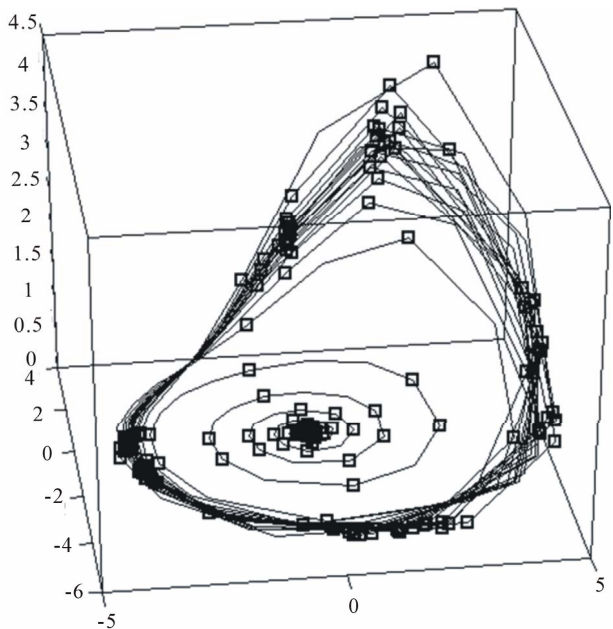

**Figure 6. Attractor with markers.**

The presence of the symmetry of the system determined the structure of the form (3). Parameter identification system using the method of least squares, gives the following result:

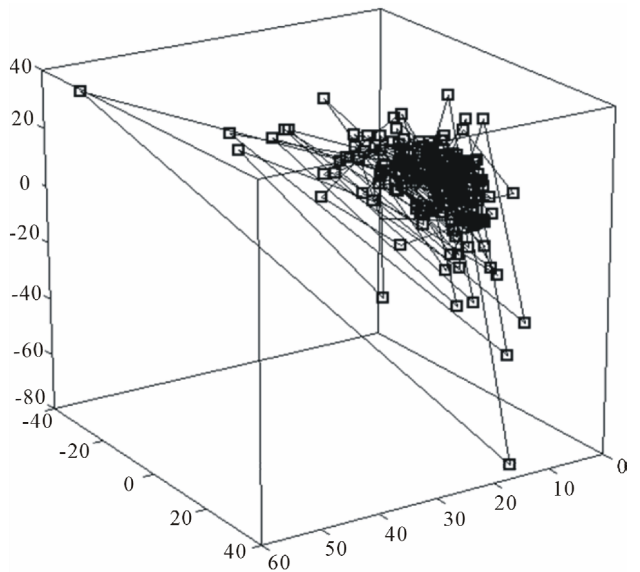

**Figure 8. Labeled attractor.**

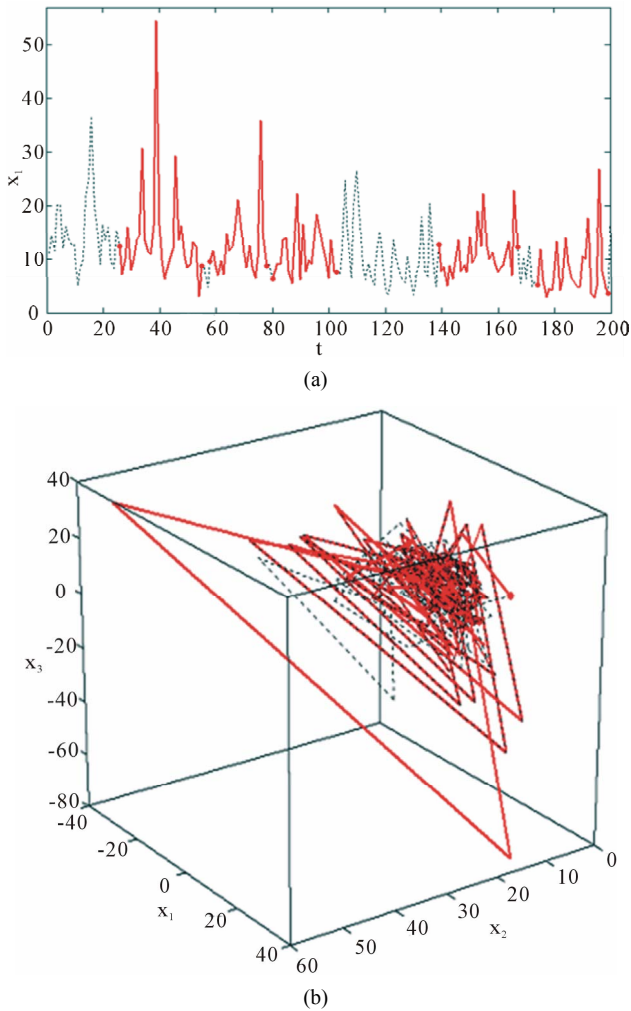

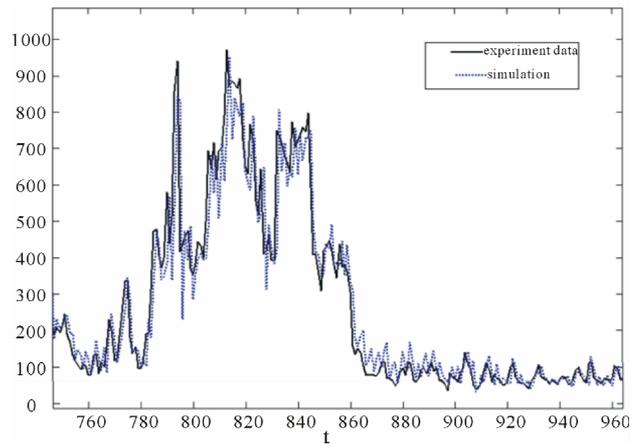

Figure 10. Comparison of dynamics.

sion can be close enough to the identity for small perturbations.

## 6. Conclusion

The use of symmetric properties of evolution equations can be successfully used for the reconstruction of the system. Developing a constructive method for structural identification systems with chaotic dynamics based on geometric analysis of the experimentally obtained phase portrait of the system.

Figure 9. Decision-winner at baseline (a) and the attractor (b).

$$A = \begin{bmatrix} 0.9413 & -0.1805 & 0.1164 & -0.0295 \\ -0.0545 & 0.8226 & 0.1622 & 0.1056 \\ 0.0014 & -0.0105 & -0.4455 & 0.8471 \\ -0.0062 & 0.0341 & -0.8860 & -0.5404 \end{bmatrix},$$

$$\Psi_0 = \begin{bmatrix} 0.0399 \\ 0.0463 \\ -0.4848 \\ -0.1851 \end{bmatrix} \left( \exp\left(t^{0.0001}\right) \sin\left(t^{0.4}\right) \right),$$

$$C = 10^4 \begin{bmatrix} 2.1037 & -0.0124 & 0.1202 & -0.0302 \end{bmatrix}.$$

**Figure 10** shows a comparison of the dynamics of the original system and the reconstructed model.

Note that in determining the conversion makes sense to consider the conservation of structural stability for flows. Based on the equivalence of all the disturbances in model assumed that the local topological equivalence should preserve the structural stability, and the conver-